\documentclass[runningheads]{llncs}

\usepackage{algorithm}
\usepackage{algorithmicx}
\usepackage{algpseudocode}
\usepackage{verbatim}
\usepackage{blkarray}
\usepackage{amsmath}

\usepackage{amsfonts}
\usepackage{tikz-cd}
\usepackage{xcolor}

\renewcommand{\O}{\mathcal{O}}
\renewcommand{\u}{{\bf u}}
\newcommand{\x}{{\bf x}}
\newcommand{\p}{{\bf p}}
\newcommand{\mydef}[1]{{\textit{#1}}}

\usepackage[T1]{fontenc}
\usepackage{graphicx}
\begin{document}
\title{Monodromy Coordinates}
%
%
\author{Taylor Brysiewicz\inst{1}\orcidID{0000-0003-4272-5934} }
\authorrunning{T. Brysiewicz}
%
\institute{University of Western Ontario, London ON N6G 2V4, Canada 
\email{tbrysiew@uwo.ca}}
\maketitle              
\begin{abstract}
We introduce the concept of monodromy coordinates for representing solutions to large polynomial systems. Representing solutions this way provides a time-memory trade-off in a monodromy solving algorithm. We describe an algorithm, which interpolates the usual monodromy solving algorithm, for computing such a representation and analyze its space and time complexity. 
\end{abstract}
\section{Introduction} 
\begin{figure}[!htpb]
\includegraphics[scale=0.36]{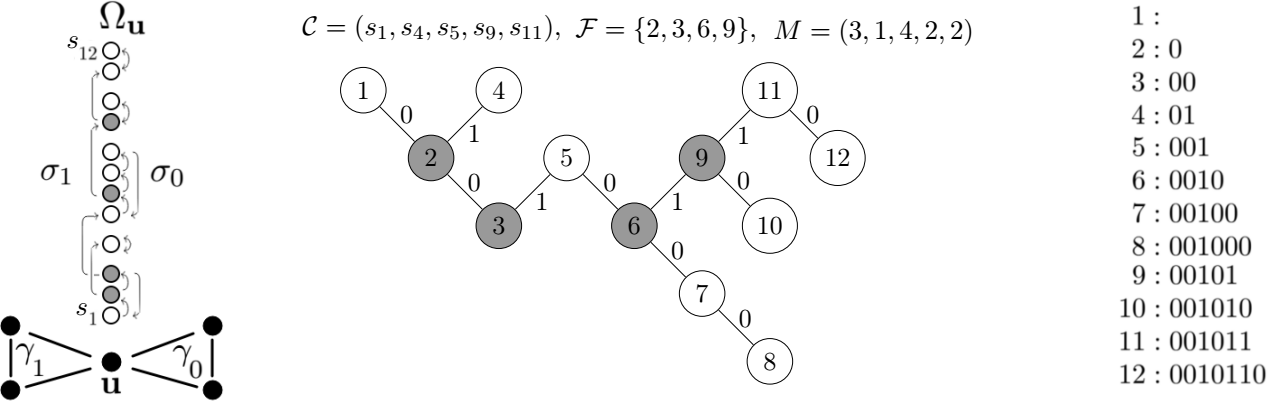}
\caption{A depiction of monodromy (left) a monodromy tree of type \textrm{IV} (center) and monodromy coordinates (right), for a system with twelve solutions.}
\label{fig:firstpicture}
\end{figure}

The bottleneck for numerically solving a polynomial system can lie in the space-complexity of an algorithm rather than its time-complexity: storing $d$ points in $\mathbb{C}^n$ requires $\O(dn)$ bits.  For instance, the solution set to the system in Example~\ref{ex:bigsystem} requires around $100$ gigabytes to store. We propose the alternative representation of a \textit{monodromy tree} which encodes the \textit{monodromy coordinates} of the  solutions (see Figure \ref{fig:firstpicture}). A monodromy tree describes how to find  solutions using a  \textit{monodromy solving algorithm} (Algorithm~\ref{alg:monodromy_solve}) and provides an iterator (Algorithm~\ref{alg:next}) for the solution set. The expected space complexity of storing a monodromy tree is $\O(n\ln(d))$.  Using a monodromy tree, the solutions in Example~\ref{ex:bigsystem} can be represented using less than $15$ kilobytes (see Example \ref{ex:bigsystem}).

We analyze our work under the three {\textit{strong monodromy assumptions}} (see Section \ref{sec:Monodromy}). These place our analysis within the realm of probabilistic group theory. Hence, our complexity estimates refer to \textit{expected} behaviour. 
The memory reduction achieved through monodromy coordinates is balanced by the $d$ queries to \textit{monodromy oracles}  required to \textit{unpack} the solutions (see Theorem~\ref{thm:storage}), which is comparable to the $\O(d)$ queries for solving the system in the first place. However, since a monodromy tree gives an  iterator for the solution set (Algorithm~\ref{alg:next}), accumulating functions of the solutions (e.g., the number of real/complex solutions) can be achieved without ever holding all solutions in memory.

Crucially, a monodromy tree can be computed using an augmented monodromy solving algorithm (Algorithm \ref{alg:monodromy_solve_MTE}) which has the same space complexity  as storing a monodromy tree ($\O(n\ln(d))$), but comes at the cost of a time complexity increase   to $\O(d(d/\ln(d)))$. Optionally, one may interpolate between Algorithms \ref{alg:monodromy_solve} and \ref{alg:monodromy_solve_MTE} based on a parameter $1\leq \alpha \leq d$. This gives an algorithm with space complexity $\O(n(d/\alpha))$ and time complexity $\O(d\alpha)$. As a consequence, our modified monodromy solving algorithm induces a time-memory trade-off, quantified by a factor of~$\alpha$.

\section{Monodromy Background}
\label{sec:Monodromy}
Fix a  system in variables $\mydef{\x}=(x_1,\ldots,x_n)$ and parameters $\mydef{\p}~=~(p_1,\ldots,p_k)$:
\[
F(\x;\p) := ({f_1(\x;\p)},\ldots,{f_N(\x;\p)}) \subset \mathbb{C}[p_1,\ldots,p_k][x_1,\ldots,x_n] = \mathbb{C}[\p][\x]
\]
such that the incidence variety $
\mydef{X} = \{(\x,\p) \mid F(\x;\p)=0\}\subseteq \mathbb{C}^n \times \mathbb{C}^k$ 
is irreducible and the image of the projection ${\pi}:X \to \mathbb{C}^k$ has the same dimension as $X$. Let ${Y}:=\overline{\pi(X)}$ be the Zariski closure of this image. Then $\pi:X \to Y$ is a branched cover of degree $\mydef{d}$ for some $d \in \mathbb{N}$ and there exists a Zariski open subset ${\mathcal U} \subseteq Y$ of \mydef{regular values} for which $\pi|_{\pi^{-1}(\mathcal U)}$ is a $d$-to-one covering space. 

 Given a regular value $\u \in \mathcal U$, any loop ${\gamma}:[0,1] \to \mathcal U$ based at $\u$ lifts to $d$ paths $\pi^{-1}(\gamma([0,1]))$. These paths induce a permutation ${\sigma_\gamma}$ on the fibre  $\pi^{-1}(\u) = {\Omega_{\u}}$ mapping start points to end points. Write ${G_{\u}}$ for the image of the homomorphism $\gamma \mapsto \sigma_{\gamma}$ from the fundamental group of $\mathcal U$ to the symmetric group ${\mathfrak S_{\Omega_\u}}$. The loop $\gamma$ is a \mydef{(monodromy) loop}, the permutation $\sigma_{\gamma}$ is a \mydef{(monodromy) permutation}, and the group $G_{\u}$ is the \mydef{monodromy group} of $\pi$ based at $\u$.

We analyze our work under the \mydef{strong monodromy assumptions}:\begin{itemize}
\item[] \textbf{Full Symmetry:} $G_\u = \mathfrak S_{\Omega_\u} \cong \mathfrak S_d$
\item[] \textbf{Uniform Sampling:}  We can sample elements of $G_\u$ uniformly at random
\item[] \textbf{Evaluation Oracle:}  Given $\sigma \in G_\u$ and $s \in \Omega_\u$ we can evaluate $\sigma(s)$
\end{itemize}  
These assumptions are reasonable. (Full-Symmetry) Many  branched covers have full-symmetric monodromy groups \cite{Harris,Yahl}. (Uniform Sampling) Given a distribution on $\mathfrak S_{d}$ with no zero probabilities, taking convolutions gives distributions converging to the uniform distribution \cite[Ch 4.~Thm 3.]{Diaconis}. (Evaluation Oracle) The numerical method of \mydef{homotopy continuation} is incredibly reliable and functions as a practical evaluation oracle for monodromy. For details about homotopy continuation and numerical algebraic geometry, see \cite{NumericalNonlinearAlgebra,IntroNumericalAlgGeom:2005}. 

The strong monodromy assumptions are often made when analyzing monodromy algorithms because they place the  analysis of these algorithms directly within the realm of \textit{probabilistic group theory}. To that end, we recall some facts regarding random permutations.

The first result, due to Dixon \cite{Dixon}, states that the probability ${p_d}$ that two random permutations $\sigma_0,\sigma_1 \in \mathfrak{S}_{d}$ generate a transitive subgroup is $1-d^{-1}+\O(d^{-2})$. For reference, $p_{1000}>0.998998$. The second fact, is that the expected number of $i$-cycles in a random permutation $\sigma \in \mathfrak S_d$ is $\frac{1}{i}$. Consequently, the expected number of cycles in a permutation is 
 given by the \mydef{$d$-th Harmonic number}
$
{H_d}~=~\sum_{i=1}^d \frac{1}{i} = \frac{1}{1} + \frac{1}{2} +\cdots+\frac{1}{d} \approx \ln(d).
$  Finally, the expected largest cycle in a random permutation is asymptotically $\lambda \cdot d$ where ${\lambda} \approx 0.6243$ is the \mydef{Golomb-Dickman constant}~\cite{Golomb}.

\section{Monodromy Solving}
Monodromy plays a central role in the algorithms and applications of modern numerical algebraic geometry \cite{Monodromy:2019,NumericalMonodromy,Yahl}. A \mydef{monodromy solving algorithm} computes the $d$ solutions $\Omega_{\u}$ to a system $F(\x;\u)=0$ from three  inputs: (a) a \mydef{seed solution} $s_1 \in \Omega_{\u}$, (b) evaluation oracles ${\sigma}$ for $D$ monodromy permutations $\sigma_0,\ldots,\sigma_{D-1} \in \mathfrak S_{\Omega_{\u}}$, and (c) a \mydef{\texttt{stopping criterion}}.

  A monodromy solving algorithm succeeds as long as $\langle \sigma_0,\ldots,\sigma_{D-1}\rangle$ is a transitive permutation group. This is incredibly likely for $D=2$ and $d\gg 0$ by Dixon's result. Hence, we consider only two random permutations ${\sigma_0},{\sigma_1} \in \mathfrak S_{\Omega_{\u}}$ and give a deterministic {monodromy solve} algorithm: Algorithm~\ref{alg:monodromy_solve}. For a more versatile framework for monodromy algorithms, see \cite{Monodromy:2019}.

\begin{algorithm}[htb]
\small
\caption{{\sf MonodromySolve}}\label{alg:monodromy_solve}
\algnotext{EndIf}
\algnotext{EndFor}
\algnotext{EndWhile}
\begin{flushleft}
\textbf{Input: }{$\bullet$ A seed solution $s_1 \in \Omega_\u$ \quad \quad  $\bullet$ Evaluation oracles $\sigma$ for $\sigma_{0},\sigma_{1} \in \mathfrak S_{\Omega_\u}$ which generate a transitive permutation action \quad \quad $\bullet$ A \texttt{stopping criterion}}  \\
\textbf{Output: }{The $d$ points of $\Omega_\u$} 
\end{flushleft}
\begin{algorithmic}[1]
\vspace{-10pt}
\State \texttt{initialize} $S = \{s_1\}$ and $j=1$ \Comment{Initialize the set of found solutions} 
\While{$\texttt{stopping criterion} = \texttt{false}$}
\If{$q:=\sigma_0(\texttt{last}(S))$ is new} \Comment{Check if a new point is found via $\sigma_0$ }
    \State $S \leftarrow q$
\Else{
    \While{$q=\sigma_1(s_j)$ is not new} \Comment{If not, try finding new points via $\sigma_1$}
    	\State $j=j+1$
    \EndWhile
    \State $S \leftarrow q$;\quad  $j = j+1$ \Comment{Such $j$ are called founder indices}
}
\EndIf
\EndWhile
\State {\hspace{15pt}\Return{$S$}}\Comment{The order of $S$ gives the monodromy ordering of $\Omega_{\u}$}
\end{algorithmic}
\end{algorithm}

Algorithm \ref{alg:monodromy_solve} greedily discovers new solutions by applying $\sigma_0$ to known solutions, until the  known cycles of $\sigma_0$ have been saturated. It then applies $\sigma_1$ to the known solutions, in the order they have been found, until a new cycle of $\sigma_0$ is discovered, and repeats until the stopping criterion is met. 
For now, we take the stopping criterion in Algorithm \ref{alg:monodromy_solve} to be \textit{$|S|=d$}.

Algorithm \ref{alg:monodromy_solve} induces an ordering on $S = \Omega_{\u}$ called the \mydef{monodromy ordering} of $\Omega_{\u}$ (with respect to $\sigma_0,\sigma_1,s_1$). We write ${\mathcal C}:=({c_1},\ldots,{c_r})$ for the minimal elements of each cycle of $\sigma_0$ with respect to this order, called \mydef{initial cycle solutions}, and ${M}:=({m_1},\ldots,{m_r})$ for the corresponding cycle sizes. In particular, ${r}$ is the number of cycles of $\sigma_0$. 
 Each cycle  is \textit{found} by some solution $s_j$ via $\sigma_1(s_j)$ in line $\texttt{8}$ in which case we call $j$ a \mydef{founder index}. We denote the set of $r-1$ founder indices by ${\mathcal F}$, whose largest element we denote by ${j^*}$.  For ${\alpha} \in \mathbb{N} \cup \{\infty\}$, we write ${\mathcal C_{\alpha}}$ for the tuple of initial cycle solutions, along with every $\alpha$-th solution in that cycle (e.g. $\mathcal C_{\infty} = \mathcal C$ and $\mathcal C_1 = \Omega_{\u}$).

\begin{example}
\label{ex:permutations}
Consider  $\Omega_{\u}=\{s_1,\ldots,s_{12}\}$ along with the permutations 
\begin{align*}
\sigma_0&=(s_1,s_2,s_3)(s_4)(s_5,s_6,s_7,s_8)(s_9,s_{10})(s_{11},s_{12}) \\ \sigma_1 &= (s_1,s_3,s_5,s_2,s_4)(s_6,s_9,s_{11})(s_7)(s_8,s_{10})(s_{12}) \in \mathfrak S_{\Omega_{\u}}
\end{align*}
\noindent shown in Figure \ref{fig:permutations}.

\begin{figure}[!htpb]
\begin{center}
\begin{tikzpicture}[node distance=2cm, scale=0.96]

\node (N) at (-1, 0) {$\sigma_0$:};
\node (1) at (0, 0) {(\,\,$s_1$};
\node[draw] (2) at (1, 0) {$s_2$};
\node[draw] (3) at (2, 0) {$s_3$};
\node (n) at (2.4, 0) {)};
\node (4) at (3, 0) {(\,\,$s_4$\,\,)};
\node (5) at (4, 0) {(\,\,$s_5$};
\node[draw] (6) at (5, 0) {$s_6$};
\node (7) at (6, 0) {{\color{white}(}$s_7$};
\node (8) at (7, 0) {{\color{white}(}$s_8$\,\,)};
\node (nn) at (7.6, 0) {(};
\node[draw] (9) at (8, 0) {$s_9$};
\node (10) at (9, 0) {{\color{white}(}$s_{10}$\,\,)};
\node (11) at (10, 0) {(\,\,$s_{11}$};
\node (12) at (11, 0) {{\color{white}(}$s_{12}$\,\,)};

\draw[->, bend right=80, shorten >=1pt] (2) edge (4);
\draw[->, bend right=80] (6) edge (9);
\draw[->, bend right=80] (9) edge (11);
\draw[->, bend right=80] (3) edge (5);
\draw[->, bend left=55, dashed] (1) edge (3);
\draw[->, bend right=70, dashed] (4) edge (1);
\draw[->, bend right=70, dashed] (5) edge (2);
\draw[->, bend right=60, dashed] (11) edge (6);
\draw[<->, bend right=60, dashed] (10) edge (8);
\draw[<-, dashed,shorten <=1pt] (7) to [out=120,in=60,loop,looseness=4.8] (7);
\draw[<-, dashed,shorten <=1pt] (12) to [out=120,in=60,loop,looseness=4.8] (12);

\end{tikzpicture}

\caption{Two permutations $\sigma_0, \sigma_1$ on twelve solutions $\Omega_{\u}$. In the middle, $\sigma_0$ is written in cycle notation, whereas $\sigma_1$ is represented via arrows. A solid arrow from a boxed founder indicates an instance of line $\texttt{6}$ in Algorithm \ref{alg:monodromy_solve} proceeding to line $\texttt{8}$.}
\label{fig:permutations}
\end{center}
\end{figure}
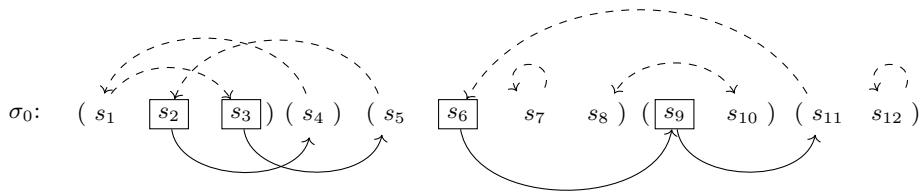

The solutions  $\{s_1,\ldots,s_{12}\}$ are already indexed with respect to their monodromy order and Algorithm \ref{alg:monodromy_solve} finds them via  $20=11_{\sigma_0}+9_{\sigma_1}$~queries.
Running Algorithm~\ref{alg:monodromy_solve} on $(\sigma_0,\sigma_1,s_1)$ gives the values $\mathcal F = \{2,3,6,9\}$, $ j^* = 9,$ $  \mathcal C~=~(s_1,s_4,s_5,s_9,s_{11}),  \mathcal C_2 = (s_1,s_3,s_4,s_5,s_7,s_9,s_{11}),$ and $M=(3,1,4,2,2).$
\end{example}

We write {$\rho(\square)$} for the number of bits used to store an object $\square$. For example, $\rho(\mathbb{Z})=64$ for moderate sized integers, and $\rho(\mathbb{C}^n) = \rho(\mathbb{R}^{2n}) = 2n\rho(\mathbb{R})$. With this notation, we describe the time and space complexity of Algorithm \ref{alg:monodromy_solve}.

\begin{proposition}
Algorithm \ref{alg:monodromy_solve} has time complexity $\O(d)$ and space complexity $\O(nd)$ as measured in oracle queries and bits respectively.
 Specifically, Algorithm~\ref{alg:monodromy_solve} queries a monodromy oracle $d+j^*-1$ times and requires $d\rho(\mathbb{C}^n) $ bits to store the solutions. 
\end{proposition}
\begin{proof}
Index $\Omega_\u$ via the monodromy ordering. Evaluation oracles are queried in lines \texttt{3} and \texttt{6}. The former is evaluated $d-1$ times, once for every solution other than $s_d$. Line \texttt{6} is evaluated for $j=1,\ldots,j^*$. This totals to $d+j^*-1$ as stated. The space required to perform the algorithm is negligibly more than the cost of storing $d$ points in $\mathbb{C}^n$ (e.g. $q$ and $j$ must be stored).
\end{proof}

\section{Monodromy Coordinates: Compression}
Each solution $s \in \Omega_{\u}$ found via Algorithm \ref{alg:monodromy_solve} occurs as a sequence of permutations, each either $\sigma_0$ or $\sigma_1$, applied to $s_1$. Identify this sequence with a word ${\omega}={\omega_1\cdots \omega_\ell}$  in $ \{0,1\}$, that is,
$
s =(\sigma_{\omega_{k}}\circ \cdots \circ \sigma_{\omega_1})(s_1) \longleftrightarrow \omega=\omega_1\cdots\omega_\ell.
$ We say that $\omega$ gives \mydef{monodromy coordinates} for $s$ with respect to $(\{\sigma_0,\sigma_1\},s_1)$. By construction, the monodromy coordinates of $s$ with respect to $(\{\sigma_0,\sigma_1\},s_1)$ describes the lexicographically smallest ($\sigma_0<\sigma_1$) path from $s_1$ to $s$ in the Cayley-graph of $\langle \sigma_0,\sigma_1 \rangle$ acting on $\Omega_{\u}$. The induced monodromy ordering of $\Omega_{\u}$ agrees with the lexicographic ordering on the monodromy coordinates.

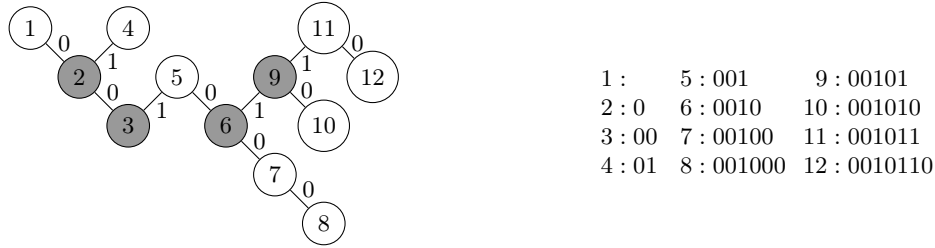
\begin{figure}[!htpb]
\minipage{0.32\textwidth}
\begin{tikzpicture}[node distance=2cm, scale=0.65]
\node[draw, circle] (1) at (0, 0) {1};
\node[draw, circle,fill=black!40!]  (2) at (1, -1) {2};
\node[draw, circle,fill=black!40!]  (3) at (2, -2) {3};
\node[draw, circle]  (4) at (2, 0) {4};
\node[draw, circle]  (5) at (3, -1) {5};
\node[draw, circle,fill=black!40!]  (6) at (4,-2) {6};
\node[draw, circle]  (7) at (5,-3) {7};
\node[draw, circle]  (8) at (6,-4) {8};
\node[draw, circle,fill=black!40!]  (9) at (5,-1) {9};
\node[draw, circle]  (10) at (6,-2) {10};
\node[draw, circle]  (11) at (6,0) {11};
\node[draw, circle] (12) at (7,-1) {12};

\draw[-] (1) edge node[pos=0.5,anchor=-135,inner sep=1pt] {0} (2);
\draw[-] (2) edge node[pos=0.5,anchor=-135,inner sep=1pt] {0}  (3);
\draw[-] (2) edge node[pos=0.5,anchor=135,inner sep=1pt] {1}  (4);
\draw[-] (3) edge node[pos=0.5,anchor=135,inner sep=1pt] {1}  (5);
\draw[-] (5) edge node[pos=0.5,anchor=-135,inner sep=1pt] {0}  (6);
\draw[-] (6) edge node[pos=0.5,anchor=-135,inner sep=1pt] {0}  (7);
\draw[-] (6) edge node[pos=0.5,anchor=135,inner sep=1pt] {1}  (9);
\draw[-] (7) edge node[pos=0.5,anchor=-135,inner sep=1pt] {0}  (8);
\draw[-] (9) edge node[pos=0.5,anchor=135,inner sep=1pt] {1}  (11);
\draw[-] (9) edge node[pos=0.5,anchor=-135,inner sep=1pt] {0}  (10);
\draw[-] (11) edge node[pos=0.5,anchor=-135,inner sep=1pt] {0}  (12);
\end{tikzpicture}
\endminipage \hspace{1.5in}
\minipage{0.32\textwidth}
\begin{tabular}{lll}
$1:$ \hspace{0.2in} & $5:001$ \hspace{0.2in} & $\,\,\,9:00101$ \\
$2:0$ & $6:0010$ & $10:001010$ \\
$3:00$ & $7:00100$ & $11:001011$ \\
$4:01$ & $8:001000$ & $12:0010110$
\end{tabular}
\endminipage
\caption{(Left) The monodromy tree associated to the permutations $\sigma_0,\sigma_1$ in Figure \ref{fig:permutations} with founders shaded.  (Right) The corresponding monodromy coordinates. }
\label{fig:binaryTree}
\end{figure}

Monodromy coordinates can be as long as $d$ letters and so directly storing may require quadratic memory in $d$. A much more efficient encoding is via a binary tree as depicted in Figure \ref{fig:binaryTree}, whose structure is given by the monodromy coordinates: an edge labelled $0$ from $j \to j'$ occurs if $s_{j'}$ is found via $\sigma_0(s_j)$ in line $\texttt{4}$ of Algorithm \ref{alg:monodromy_solve}, and an edge labelled $1$ occurs if $s_{j'}$ is found via $\sigma_1(s_j)$ in line $\texttt{8}$. The resulting labelled rooted binary tree ${\mathcal T}$ is the \mydef{monodromy tree} of $\Omega_{\u}$, with respect to $(\{\sigma_0,\sigma_1\},s)$. See Figure \ref{fig:binaryTree} for the  monodromy tree/coordinates of Example \ref{ex:permutations}.

We consider five representations of a monodromy tree: 
\begin{equation*}\label{eq:representations} \underbrace{(s_1,\sigma)}_{\textbf{{\textrm{I}: single seed}}}\hspace{-8pt} \hookrightarrow\hspace{-5pt} \underbrace{(s_1,\sigma,\mathcal F)}_{\textbf{{\textrm{II}: minimal}}}\hspace{-5pt}\hookrightarrow\hspace{-8pt}  \underbrace{(s_1,\sigma,\mathcal F,M)}_{\textbf{{\textrm{III}: with cycle sizes}}}\hspace{-5pt}\hookrightarrow\hspace{-8pt}  \underbrace{(\mathcal C,\sigma,\mathcal F, M)}_{{\substack{\textbf{\textrm{IV}: with initial cycle}\\\textbf{representatives}}}}\hspace{-5pt}\hookrightarrow\hspace{-5pt}  \underbrace{(\mathcal C_{\alpha},\sigma,\mathcal F, M)}_{{\substack{\textbf{\textrm{V}: with cycle}\\\textbf{representatives}}}}
\end{equation*}
The cost of unpacking the solution set $\Omega_{\u}$ from a monodromy tree $\mathcal T$  differs depending on how $\mathcal T$ is represented. Indeed, the \textit{single seed} representation of type \textrm{I} is merely the input of Algorithm \ref{alg:monodromy_solve} and so one recovers $\Omega_{\u}$ by running Algorithm \ref{alg:monodromy_solve}. Given a type \textrm{II} representation of $\mathcal T$, the solution set $\Omega_{\u}$ can be recovered by executing Algorithm~\ref{alg:monodromy_solve} subject to the minor, efficiency-boosting, change: replace lines \texttt{6} and \texttt{7} with $q = \sigma_1(s_{\mathcal F_j})$.
 Thus, the number of oracle queries  to recover all solutions from a type \textrm{II} representation is $
d+r-2 = (r-1)+ (d-r) + (r-1).$
If cycle sizes are included (type \textrm{III}) then   $r-1$ evaluations of $\sigma_0$ can be saved. When initial cycle solutions are stored (type \textrm{IV}), evaluations of $\sigma_1$ are entirely unnecessary, saving $r-1$ additional oracle queries. When $\lceil m_i/\alpha \rceil$-many solutions in each cycle are stored (type \textrm{V}), they need not be found via oracle calls.
\begin{theorem}
\label{thm:storage}
The memory cost, more than $\rho(\sigma)$, to store a monodromy tree is 
\begin{align*}
\textrm{I}:\,\,&  \rho(\mathbb{C}^n) \xrightarrow{\mathbb{E}} \rho(\mathbb{C}^n)\xrightarrow{\O}\O(n)  \\ 
\textrm{II}:\,\,& \rho(\mathbb{C}^n)+(r-1)\rho(\mathbb{N})\xrightarrow{\mathbb{E}}   \rho(\mathbb{C}^n)+(H_d-1) \rho(\mathbb{N})\xrightarrow{\O}\O(\max(n,\ln(d)))  \\ 
\textrm{III}:\,\,& \rho(\mathbb{C}^n)+ (2r-1) \rho(\mathbb{N})  \xrightarrow{\mathbb{E}}\rho(\mathbb{C}^n)+ (2H_d-1)\rho(\mathbb{N})  \xrightarrow{\O} \O(\max(n,\ln(d)))  \\
\textrm{IV}:\,\,&  r\rho(\mathbb{C}^n)+ (2r-1) \rho(\mathbb{N})  \xrightarrow{\mathbb{E}} H_d\rho(\mathbb{C}^n)+ (2H_d-1) \rho(\mathbb{N})\xrightarrow{\O} \O(n\ln(d))  \\
\textrm{V}:\,\,&  {\footnotesize{\left(\sum_{i=1}^r\lceil \frac{m_i}{\alpha} \rceil\right)\rho(\mathbb{C}^n)+ (2r-1) \rho(\mathbb{N}) \xrightarrow{\mathbb{E}} < ((d/\alpha)+H_d)\rho(\mathbb{C}^n)+(2H_d-1)\rho(\mathbb{N})}} \\ & \hspace{147pt} \xrightarrow{\O}\O(\max(n(d/\alpha),n\ln(d))) \\  
\end{align*}

\vspace{-15pt}

\noindent The time cost of unpacking all solutions from a monodromy tree is $\O(d)$:

\vspace{-10pt}

\[
\textrm{I}:\,\, d+j^*-1 \quad 
\textrm{II}:\,\, d+r-2 \quad 
\textrm{III}:\,\,d-1\quad 
\textrm{IV}:\,\,d-r \quad 
\textrm{V}:\,\, d-\left(\sum_{i=1}^r \lceil \frac{m_i}{\alpha} \rceil\right)
\]

\vspace{-5pt}
\end{theorem}
\begin{proof}The memory cost directly follows from the descriptions of the monodromy tree representations. The expectation of the memory applies (a) expectation on the number of cycles, $\mathbb{E}(r) = H_d \approx \ln(d)$, along with (b) the observation that the number of points stored in type $\textrm{V}$ is at most $(d/\alpha) +r$. The asymptotics are immediate.  The time cost (in oracle queries) summarizes the discussion prior to the theorem statement. Again, the expected asymptotics are immediate.
\end{proof}

\begin{remark}\label{rem:NAGapplication}
For monodromy solving in practice we use $
\rho(\mathbb{C}^n) = 2n\rho(\mathbb{R})= 128\cdot n,$ $\rho(\mathbb{N}) = 64$, and $\rho(\sigma) = 5\rho(\mathbb{C}^k) = 640k$, by storing a real number as a $64$-bit float and an integer at most $d$ as a $64$-bit integer. For extremely large $d$, one would have to use larger integer types. Note that  two monodromy loops can be represented by five parameter values via a bow-tie configuration (see Figure~\ref{fig:firstpicture}).
For $\alpha < d/r$ the memory storage for type $\textrm{V}$ in Theorem \ref{thm:storage} is dominated by $128 n\frac{d}{\alpha}$.
\end{remark}

\begin{example}\label{ex:bigsystem} There are 
 $d=666841088$ quadric surfaces tangent to $9$ given quadric surfaces in $\mathbb{C}^3$ (see \cite{TangentQuadrics}). Finding these tangent quadrics amounts to solving a polynomial system in $10$ variables, over $90$ parameters. Under the practical storage assumptions of Remark \ref{rem:NAGapplication}, we specialize Theorem \ref{thm:storage} for various monodromy tree representations. We use the approximation $H_d \approx \ln(d)$ and the empirical observation that $\mathbb{E}(j^*) \approx \frac{d+1}{2}$. Note that $d/\ln(d) \approx 32820112$. To encode a monodromy tree of type $\textrm{I}$ costs $58,880$ bits and requires around one-billion oracle calls to unpack. For types $\textrm{II}-\textrm{IV}$, the number of oracle calls is approximately $d$ and the approximate number of bits to store are $60,116, 61,416,$ and $112,150$ respective. Type $\textrm{V}$ with $\alpha = 100000,10,$ and $1$ respectively require approximately $8$-million, $85$-billion, and $853$-billion bits to store. To unpack these three trees of type $\textrm{V}$ one needs approximately $d$,$\frac{9}{10}d$, and $\frac{1}{2}d$ oracle queries respectively.
\end{example}

\section{Monodromy Solving Using Monodromy Coordinates}

 A key advantage of monodromy coordinates lies in the fact that the monodromy solve algorithm can be performed in such a way to produce a monodromy tree representation of type $\textrm{IV}$ or $\textrm{V}$. During this procedure, the algorithm ``knows'' when new solutions have been found and can thus accumulate any function ${\mathcal G(\u)} := \sum_{s \in \Omega_{\u}} g(s)$ during the process, avoiding the need to store all solutions simultaneously in memory. Examples of desirable \mydef{accumulators} include 
\begin{equation}
\label{eq:accumulators}
\textbf{\#}\mathbb{C}:  {\mathcal G_{\mathbb{C}}} = \sum_{s \in \Omega_{\u}} 1,  \quad \textbf{\#}\mathbb{R}:  {\mathcal G_{\mathbb{R}}} = \sum_{s \in \Omega_{\u} \cap \mathbb{R}^n} 1, \quad  \textbf{Trace: } {\mathcal G_{\Sigma(\textbf{v})}} = \sum_{s \in \Omega_{\textbf{v}}} s.
\end{equation}
The trace accumulator is particularly desirable when the polynomial system corresponds to the computation of a \textit{witness set}. In this case, the trace accumulator can be used as a stopping criterion via the \textit{trace test} \cite{TraceTest:2018}. Similar trace  accumulators can be designed for other trace tests, like the one in \cite{SparseTraceTest}.

A monodromy tree of type $\textrm{IV}$ gives an \textit{iterator} for the (monodromy) ordered list $\Omega_{\u}=(s_1,\ldots,s_d)$, represented by a seed solution and a  \mydef{\texttt{next}} function.

\begin{algorithm}[htb]
\small
\caption{{\sf next}}\label{alg:next}
\algnotext{EndIf}
\algnotext{EndFor}
\algnotext{EndWhile}
\begin{flushleft}
\textbf{Input: }{
 \hspace{1pt} $\bullet$ A monodromy tree $\mathcal T = (\mathcal C, \sigma, \mathcal F, M)$ of type $\textrm{V}$ \\ \quad \quad \quad\hspace{8.5pt} $ \bullet$ $j$\\ \quad \quad \quad \,\,\, \hspace{1pt}  $\bullet$ $ sj$ representing $s_j \in \Omega_{\u}$ other than $s_d$
}  \\
\textbf{Output: }{$s_{j+1}$} \\
\end{flushleft}
\begin{algorithmic}[1]
\vspace{-10pt}
\If{$j$ is a partial sum $\sum_{i=1}^{k-1} M_i$} \,\Return{$\mathcal C_{k}$} \textbf{else} \textbf{return} $\sigma_0(s_j)$
\EndIf
\end{algorithmic}
\end{algorithm}

 Algorithm \ref{alg:next} uses at most one oracle call and succeeds even when $\mathcal T$ is incomplete: for
${\mathcal T^{(k)}} := ({\mathcal C^{(k)}},\sigma,{\mathcal F^{(k-1)}},{M^{(k)}})$ where $\mathcal C^{(k)}$ is the truncation of $\mathcal C$ to those solutions only in the first $k$ cycles, $\mathcal F^{(k-1)} = (\mathcal F_i)_{i=1}^{k-1}$, $M^{(k)}=(M_i)_{i=1}^{k}$,
 and $j$ less than the sum of $M^{(k)}$, the output is correct.
 
Similarly, applied to a (possibly incomplete) monodromy tree  $\mathcal T^{(k)}$ and solution $s ~\in~\Omega_{\u}$, the procedure ${\texttt{inNewCycle}}$ determines whether $s$ belongs to a cycle of $\sigma_0$  represented by an element of $\mathcal C^{(k)}$.
We point out that  Algorithm \ref{alg:inNewCycle} uses at most $\min(\max(M),\alpha)-1$ oracle queries.

\begin{algorithm}[!h]
\small
\caption{{\sf inNewCycle}}\label{alg:inNewCycle}
\algnotext{EndIf}
\algnotext{EndFor}
\algnotext{EndWhile}
\begin{flushleft}
\textbf{Input: }{
$\bullet$ An incomplete monodromy tree $\mathcal T^{(k)} = ({\mathcal C^{(k)}},\sigma,{\mathcal F^{(k-1)}},{M^{(k)}})$ of type $\textrm{IV}$ \\  \quad \quad \quad \,\,\, $\bullet$ $ s \in \Omega_{\u}$
}  \\
\textbf{Output: }{\texttt{false} if $s$ belongs to the same cycle of $\sigma_0$ as some $c \in \mathcal C^{(k)}$, \texttt{true} otherwise} \\ \vspace{-5pt}
\end{flushleft}
\begin{algorithmic}[1]
\vspace{0pt}
\If{$s \in \mathcal C^{(k)}$} \Return \texttt{false} \EndIf
\For{$i$ from $1$ to $\textrm{min}(\textrm{max}(M^{(k)}),\alpha)-1$}
\State $s=\sigma_0(s)$
\If {$s \in \mathcal C^{(k)}$} \Return \texttt{false} \EndIf
\EndFor
\State \Return{\texttt{true}}
\end{algorithmic}
\end{algorithm}

\begin{algorithm}[!htb]
\small
\caption{{\sf MonodromySolve - Monodromy Coordinates}}\label{alg:monodromy_solve_MTE}
\algnotext{EndIf}
\algnotext{EndFor}
\algnotext{EndWhile}
\begin{flushleft}
\textbf{Input: }{\hspace{9pt}$\bullet$ A seed solution $s_1 \in \Omega_\u$ \\ \hspace{38pt} $\bullet$ Evaluation oracles $\sigma$ for $\sigma_{0},\sigma_{1} \in \mathfrak S_{\Omega_\u}$\\ \quad \quad \quad\quad\quad  (where $\{\sigma_0,\sigma_1\}$ generate a transitive action) \\ \hspace{38pt} $\bullet$ $\alpha \in \mathbb{N}$  \\\hspace{38pt} $\bullet$ A \texttt{stopping criterion}}  \\ \hspace{38pt} $\bullet$ (Optional:) any function $\mathcal G=\sum_{s \in \Omega_{\u}}g(s)$ to accumulate\\
\textbf{Output: }{A monodromy tree representation $\mathcal T = (\mathcal C, \sigma, \mathcal F, M)$ of type $\textrm{V}$} \\ \vspace{-5pt}
\end{flushleft}
\begin{algorithmic}[1]
\State $\mathcal I=\{s_1\}$, $\mathcal J=\{\}$, $M=\{\}$ \Comment{initial cycle solutions $\mathcal I$,  additional solutions $\mathcal J$, and cycle sizes $M$}
\State $m=1$, $j=1$\Comment{current cycle size $m$ and potential founder index $j$}
\State  $sj=s_1$, $ s=s_1$ \Comment{potential founder $sj$ representing $s_j$ and a placeholder $s$}
\State  $\mathcal F=\{\}$, $a=0$ \Comment{founder indices $\mathcal F$, and a storage counter $a$} 
\While{\texttt{stopping criterion} is not met}
\State \hspace{-14pt}$\star_0$\hspace{4pt}Replace $s = \sigma_0(s)$
\If{$s\neq \texttt{last}(\mathcal I)$} \Comment{Check if a new point is found via $\sigma_0$}
    \State $m=m+1, a=a+1$ \Comment{If so, increase current cycle size}
    \State $\mathcal G = \mathcal G + g(s)$ \Comment{Accumulate}
    \If{$a=\alpha$} $\mathcal J\leftarrow s; a=0$ \Comment{Store $s$ and reset the counter $a$} \EndIf
\Else{\Comment{Otherwise, a cycle is complete}
	 \State \texttt{foundNewCycle} $= \texttt{false}$ \Comment{and we must find a new one}
	\State $M \leftarrow m$, $m=1$, $a=0$ \Comment{Store cycle size in $M$ and reset counters $m$ and $a$}
	\While{\texttt{foundNewCycle} $= \texttt{false}$}
	\State \hspace{-42pt}$\star_1$\hspace{29pt} Compute $t_* := \sigma_1(sj)$ and set  $s = \texttt{copy}(t_*)$ \Comment{$sj$ is a founder if $t_*$ is new}
	\If{\hspace{-52.5pt}$\star_0$\hspace{45pt}\texttt{inNewCycle}$(\mathcal I \cup \mathcal J,\sigma,\mathcal F,M,s)$} \Comment{If this is so}
		\State \texttt{foundNewCycle} $= \texttt{true}$
		\State $\mathcal I \leftarrow t_*$, $\mathcal F \leftarrow j$ and $s = t_*$ 
		\State $\mathcal G = \mathcal G +g(s)$ \Comment{Accumulate}
	\Else
	\State \hspace{-56pt}$\star_0$\hspace{46pt}$sj=\texttt{next}(\mathcal I,\sigma,\mathcal F, M,j,sj)$ and $j=j+1$ 
	 
	\EndIf
	
	\EndWhile
}
\EndIf
\EndWhile
\State {\hspace{15pt}\Return{$\mathcal T = (\mathcal I, \sigma,\mathcal F, M)$}} \Comment{Replace $\mathcal I$ with $\mathcal C_{\alpha} = \mathcal I\cup\mathcal J$ if desired}
\end{algorithmic}
\end{algorithm}

\begin{theorem}
Algorithm \ref{alg:monodromy_solve_MTE} correctly evaluates any optional accumulators $\mathcal G$. During which, no more than $
d-1+j^*\cdot (\min(\max(M),\alpha)+1)$
oracle queries are performed. The memory cost is approximately as much as a monodromy tree of type $\textrm{V}$. In expectation,  Algorithm \ref{alg:monodromy_solve_MTE} has complexity 
\[
\text{Space: }\,\,\O(\max(n(d/\alpha),n\ln(d)))\quad \quad \quad  \text{Time: }\,\,\O(d \min(\lambda d,\alpha)).
\]
\end{theorem} 
\begin{proof}
Algorithm \ref{alg:monodromy_solve_MTE} follows Algorithm \ref{alg:monodromy_solve} with the modifications required to ascertain whether the solution obtained in line  $\texttt{6}$ of Algorithm \ref{alg:monodromy_solve} is ``new'' (using Algorithm \ref{alg:inNewCycle}) and to iterate $s_j$ (using Algorithm \ref{alg:next}). Otherwise, the only other additional operations and storage requirements come from various book-keeping steps and accumulating the functions $\mathcal G$.

The lines $\texttt{6,15,16,21}$ contribute to the oracle queries $\sigma_0$ or $\sigma_1$, as indicated by $\star_0$ and $\star_1$ decorations, respectively. Line $\texttt{6}$ queries $\sigma_0$ exactly $d-1$ times. Lines \texttt{15,16,} and \texttt{22} are performed $j^*$ times.  Only line \texttt{16} uses more than one query: it uses at most $\textrm{min}(\max(M),\alpha)-1$, totalling 
\begin{align*}
\underbrace{d\hspace{-2pt}-\hspace{-2pt}1}_{\sigma_0 \text{ on }\texttt{6}}+\hspace{-2pt}\underbrace{j^*}_{\sigma_1\text{ on }\texttt{15}}\hspace{-8pt} +  \underbrace{j^* (\min(\max(M),\alpha)-1)}_{\sigma_0 \text{ on }\texttt{16}}+\hspace{-8pt}\underbrace{j^*}_{\sigma_0 \text{ on }\texttt{21}} \hspace{-10pt}&= d\hspace{-2pt}-\hspace{-2pt}1+j^*(1\hspace{-2pt}+\hspace{-2pt}\min(\max(M),\alpha))\hspace{-2pt} 
\end{align*}
\end{proof}
In expectation, this is $\O(d\min(\lambda d,\alpha))$.
If $\alpha$ is less than $\frac{d}{\ln(d)}$, then  the expected space-time complexity estimates of Algorithm~\ref{alg:monodromy_solve_MTE} become $\O(n(d/\alpha))$ and $\O(d\alpha)$ respectively, and so Algorithm \ref{alg:monodromy_solve_MTE} gives a space-time trade-off of a factor of $\alpha$ compared to Algorithm \ref{alg:monodromy_solve}, i.e. when $\alpha=1$.

\begin{credits}
\subsubsection{\ackname}
Supported by an NSERC Discovery grant (RGPIN-2023-03551).

\end{credits}
%
%
%
%

\end{document}